\documentclass{amsart}
\usepackage{amsmath}
\usepackage{amssymb}
\usepackage[mathscr]{eucal}
\usepackage{amsthm}
\usepackage{graphicx}
\usepackage{enumerate}
\usepackage{enumitem}
\usepackage{hyperref}
\usepackage{ulem}
\usepackage{marginnote}
\usepackage{enumitem}
 \usepackage{epsfig}

\theoremstyle{plain} 
\newtheorem{theorem}{Theorem}[section]

\newtheorem{corollary}[theorem]{Corollary}
\newtheorem{definition}[theorem]{Definition}

\newtheorem{remark}[theorem]{Remark}

\def\be{\begin{equation}}
\def\ee{\end{equation}}
\def\nn{\nonumber}
\def\bea{\begin{eqnarray}}
\def\eea{\end{eqnarray}}
\def\la{\langle}
\def\ra{\rangle}

\def\IP{\hbox{\rm I\kern -1.6pt{\rm P}}}
\def\IC{{\hbox{\rm
C\kern-.58em{\raise.53ex\hbox{$\scriptscriptstyle|$}}
    \kern-.55em{\raise.53ex\hbox{$\scriptscriptstyle|$}} }}}

\def\IN{\hbox{I\kern-.2em\hbox{N}}}
\def\IR{\mathbb{R}}
\def\ZZ{\mathbb{Z}}
\def\IT{\hbox{\rm T\kern-.38em{\raise.415ex\hbox{$\scriptstyle|$}} }}

\def\IS{{\hbox{\rm
S\kern-.58em{\raise.53ex\hbox{$\scriptscriptstyle|$}}
    \kern-.55em{\raise.53ex\hbox{$\scriptscriptstyle|$}} }}}

\def\sd{\mathbf S^{d-1}}

\def\itD{\mathcal D}
\def\itW{\mathcal W}
\def\itJ{\mathcal J}
\def\itT{\mathcal T}
\def\itS{\mathcal S}
\def\itT{\mathcal T}

\def\be{\begin{equation}}
\def\ee{\end{equation}}
\def\nn{\nonumber}
\def\bea{\begin{eqnarray}}
\def\eea{\end{eqnarray}}
\def\la{\langle}
\def\ra{\rangle}

\def\IP{\hbox{\rm I\kern -1.6pt{\rm P}}}
\def\IC{{\hbox{\rm
C\kern-.58em{\raise.53ex\hbox{$\scriptscriptstyle|$}}
    \kern-.55em{\raise.53ex\hbox{$\scriptscriptstyle|$}} }}}
\def\IN{\mathbb{N}}
\def\IR{\mathbb{R}}
\def\ZZ{\mathbb{Z}}

\def\IT{\hbox{\rm T\kern-.38em{\raise.415ex\hbox{$\scriptstyle|$}} }}

\def\IS{{\hbox{\rm
S\kern-.58em{\raise.53ex\hbox{$\scriptscriptstyle|$}}
    \kern-.55em{\raise.53ex\hbox{$\scriptscriptstyle|$}} }}}

\def\sd{\mathbf S^{d-1}}

\def\cO{{\mathcal O}}
\def\cQ{{\mathcal Q}}
\def\cV{{\mathcal V}}

\def\mp{\mu^\partial}

\newtheorem{condition}{Condition}
\newtheorem{problem}{Problem}

\begin{document}


\title{Multidimensional hyperbolic billiards}

\dedicatory{Dedicated to the memory of Kolya Chernov.\\ \vskip2mm
The theory of hyperbolic billiards having been created by Sinai,\\  Kolya did more than anyone else to found their multidimensional theory.}

\author{Domokos Sz\'asz}
\address{Budapest University of Technology and Economics,
Mathematical Institute, Egry J. u. 1, 1111 Budapest, Hungary}

\email{domaszasz@gmail.com}

\date\today

\begin{abstract}
The theory of planar hyperbolic billiards  is already quite well developed by having also achieved spectacular successes. In addition there also exists an excellent monograph by Chernov and Markarian on the topic.  In contrast, apart from a series of works culminating in Sim\'anyi's remarkable result on the ergodicity of hard ball systems and other sporadic successes, the theory of  hyperbolic billiards in dimension 3 or more is much less understood. The goal of this work is to survey the key results of their theory and highlight some central problems which deserve particular attention and efforts.
\end{abstract}
\maketitle

\section{Introduction}

Mathematical billiards appeared as early as in 1912, 1913 in the works of the couple  Ehrenfest, \cite{EE12} (the wind tree model) and of D. K\"onig and A. Sz\H ucs, \cite{KSz13} (billiards in a cube) and in 1927 in the work of G. Birkhoff, \cite{B27} (those in an oval).  Ergodic theory itself owes its birth to the desire to provide mathematical foundations to Boltzmann's celebrated  ergodic hypothesis. I briefly went over its history in my survey \cite{Sz00} and more recently in the article \cite{Sz14} written on the occasion of Sinai's Abel Prize. Therefore for historic details I recommend the interested reader to consult that freely available article. Here I will only mention some of the most relevant facts from it.

In particular,
the two most significant problems from physics motivating the  initial study of mathematical billiards were
\begin{enumerate}
\item the ergodic hypothesis and
\item the goal to understand Brownian motion from microscopic principles.
\end{enumerate}
The last decades have seen further significant challenges related to billiards which have attracted the attention of both mathematicians and physicists. Without aiming at completeness we mention quantum billiards and further on promising billiard models of heat transport that have also become attackable, cf. \cite{GG08,BNSzT15}.

Beside the mathematical explanation of these fundamental problems the rigorous approach can shed light on a variety of specific questions of physics. In subsection \ref{subsub:inter} some examples will be provided to illustrate the beneficial interaction of mathematics and physics in our topic.

From the side of mathematics the 1960's saw the birth and rapid development of the theory of {\it smooth hyperbolic dynamical systems} with Sinai being one of the leading creators of this theory. Then Sinai's 1970 paper \cite{S70} introduced dispersing billiards. They represented a new object for study. In fact, they are {\it hyperbolic dynamical systems with singularities}.  Later it also turned out that this theory was  closely related to basic models of {\it chaos theory}, like logistic maps, the H\'enon map, the Lorenz map,  \dots (see for instance \cite{Y98,ChY00}). An extra push for the study of hyperbolic billiards was the revelation that, as it was not hard to see, Hamiltonian systems of hard balls with elastic collisions were isomorphic to billiards that were expected to be hyperbolic. Consequently, it became a realistic goal to prove the ergodicity of these models that would mean the first substantial step toward understanding and establishing Boltzmann's ergodic hypothesis.


As to basic notions of the theory of non-uniformly hyperbolic dynamical systems we refer to \cite{BP07}; some questions of hyperbolic billiards were also treated in \cite{KS86}. General references on ergodic theory and  dynamical systems are for instance \cite{CFS82,KH95}. For multidimensional billiards useful references are \cite{BChSzT03} and \cite{ChM06}, too; though the latter one essentially treats planar billiards only, it is also a useful reference for our main topic.

\subsection{About the paper}
For planar hyperbolic billiards there is an excellent monograph \cite{ChM06}, while for multidimensional ones nothing similar exists. (In this paper 'multidimensional' will always mean  to have dimension three or more.)  Therefore the main goal of this paper is to provide a succinct summary of the major results obtained for them.
As it, in general, happens with mathematical theories, the technique of this theory has also been developing gradually.  Therefore we also wish to pinpoint those references where essential elements of the technique are explained in a sufficiently mature way. Finally, the not too long history of the theory of multidimensional billiards shows that it leads to a variety of new phenomena and to a lot of  challenging puzzles. Here I will also mention some of these.

In section \ref{sec:summ} we summarize some key notions and main results of the theory. In particular, in subsections \ref{subsec:objects} we introduce our central objects: billiards and Lorentz processes and in \ref{subsec:prop} we list the prime stochastic properties to be understood. Subsection \ref{subsec:planekey}, where key results of the theory of planar hyperbolic are reviewed, serves as a preparation for subsection \ref{subsec:multikey} where principal results of the multidimensional theory are discussed. As said above the difference between them is that while in the planar case the theory -- if one can claim that at all -- is relatively complete and, in addition, there also exists the excellent monograph \cite{ChM06}, at the same time nothing similar can be claimed about the multidimensional theory. Indeed, in the forthcoming sections  we go into more details about the main results and some vital problems of the multidimensional theory by treating their qualitative and quantitative theory in sections  \ref{sec:multiqual} and \ref{sec:multiquan}, respectively. Finally in the Appendix those achievements of Kolya Chernov are extracted which, in my opinion, were his most significant contributions in creating and advancing the theory of multidimensional billiards.

\section{A summary}
\label{sec:summ}

\subsection{Billiards, Lorentz process}
\label{subsec:objects}

A {\it billiard} is a dynamical system describing the motion of a point particle in
a connected, compact domain $Q \subset \IT^d = \mathbb R^d / \mathbb Z^d$. In general, the boundary $\partial Q$ of
the
domain is assumed to be piecewise $C^3$-smooth; denote its smooth pieces by $\{\partial Q_\alpha| 1 \le \alpha \le J < \infty\}$. Inside $Q$ the motion is uniform while the
reflection at the boundary $\partial Q = \cup_\alpha \partial Q_\alpha$ is elastic (by the classical rule ``the angle of incidence is
equal to the angle of reflection'' hence also called specular). As the absolute value of
the velocity is a first integral of motion, the phase space of the
billiard
flow is fixed as $M=Q\times S_{d-1}$ -- in other words, every phase point
$x$ is of the form $x=(q,v)$ with $q\in Q$ and $v\in \IR^d,\ |v|=1$.
The Liouville probability measure $\mu$ on $M$ is essentially
the product of Lebesgue measures,
i.~e. $d\mu= {\rm const.}\, dq dv$ (here the constant is $\frac{1}{{\textrm vol} Q\ {\textrm vol} S_{d-1}}$) and it is invariant under the {\it billiard flow} $\{S^t: M \to M |t \in \IR\}$.

Let $n(q)$ denote the unit normal vector, directed inwards $Q$, at the point $q \in \partial Q$ within a smooth component of the
boundary $\partial Q$ .
Throughout the sequel we restrict our attention to {\it dispersing
billiards} or sometimes more broadly to {\it semi-dispersing} ones, i.~e.  we require that for every $q\in \partial Q$ the second fundamental
form $K(q)$ of the boundary component be positive (in fact, uniformly bounded away from $0$) or non-negative, respectively. The name {\it Sinai billiard} is widely known for dispersing billiards; the usage is not coherent with some authors reducing the notion to dispersing billiards without corner points while others to planar dispersing billiards.

The boundary $\partial Q$ defines a natural cross-section for the billiard
flow. Consider namely
$$
\partial M = \{ (q,v) \ | \ q\in \partial Q, \ \la v,n(q)\ra \ge 0 \} = \partial Q \times S_{d-1}^+
$$
where $S_{d-1}^+$ stands for the hemisphere of outgoing unit velocities.
The {\it billiard map}\
 $T$ is defined as the first return
map onto $\partial M$. The smooth invariant measure for the map $T$ is denoted by $\mp$,
and we have $d\mp = {\rm const.} \, |\la v,n(q)\ra |\, dq dv$ (with  ${\rm const.} = \frac{2}{{\textrm vol} \partial Q\ {\textrm vol} S_{d-1}}$).
Semi-dispersing billiards have more or less strong hyperbolic properties and consequently under some mild conditions they possess local stable and unstable invariant manifolds: $W^s_{\rm loc}$ and $W^u_{\rm loc}$ exist almost everywhere (wrt to $d\mu$ or to $d\mp$). Since billiards are Hamiltonian systems, the maximal dimensions of $W^s_{\rm loc}$ and of $W^u_{\rm loc}$ are both $d-1$ and in good cases their dimensions are indeed equal to $d-1$ a.~e. Added to that they also enjoy pleasant regularity properties for instance those of the holonomy map, distortion estimates, etc.

\begin{figure}[h]
\includegraphics[scale=0.3]{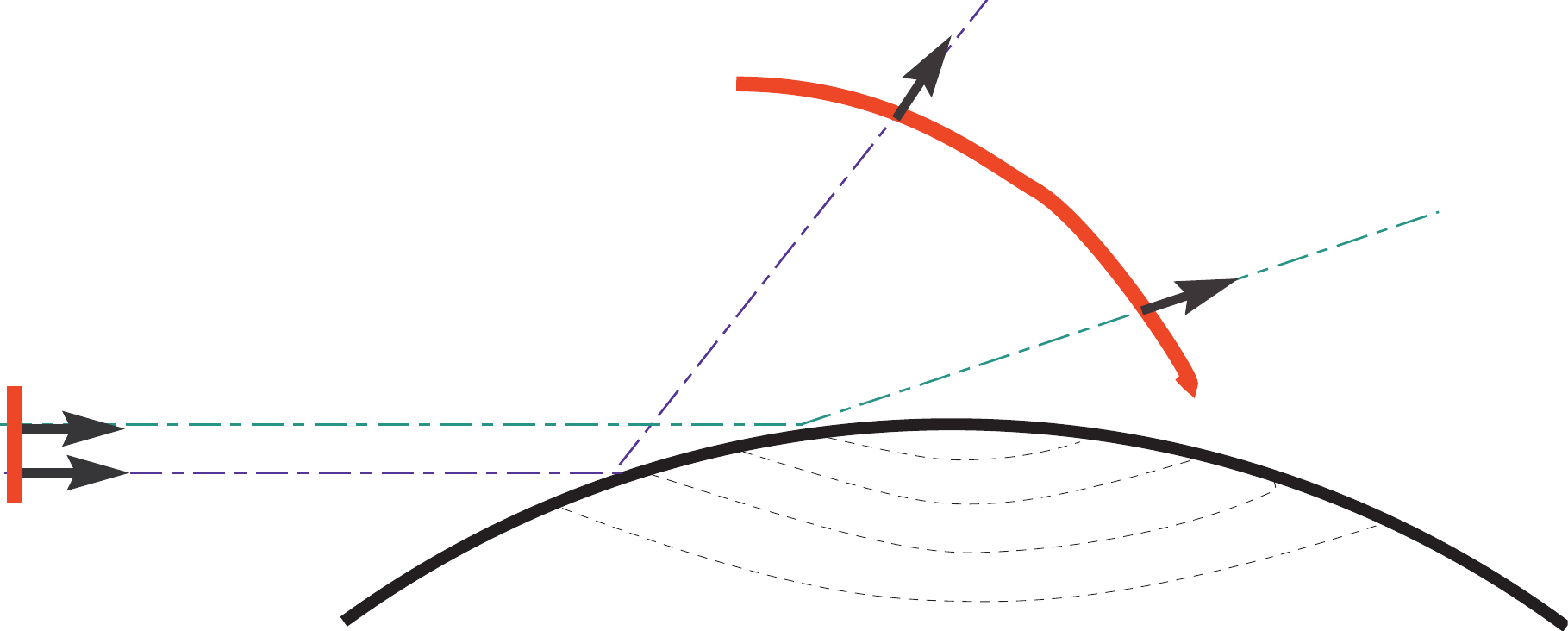}
\caption{Strong expansion}\label{fig:strong}
\end{figure}

{\it For simplicity we assume that the boundaries of the scatterers are smooth, i. e.~ there are no corner points.} As said before, semi-dispersing billiards are hyperbolic systems with singularities. Moreover, at so-called tangential singularities (cf. \ref{subsec:sing}) the derivative of the billiard ball map explodes, i. e.~ it is unbounded.  This, however, happens in  only one direction representing a strong anisotropy at collisions and  leads to additional difficulties (cf. Figures \ref{fig:strong} and  \ref{fig:weak}).

\begin{figure}[h]\label{fig:weak}
\includegraphics[scale=0.4]{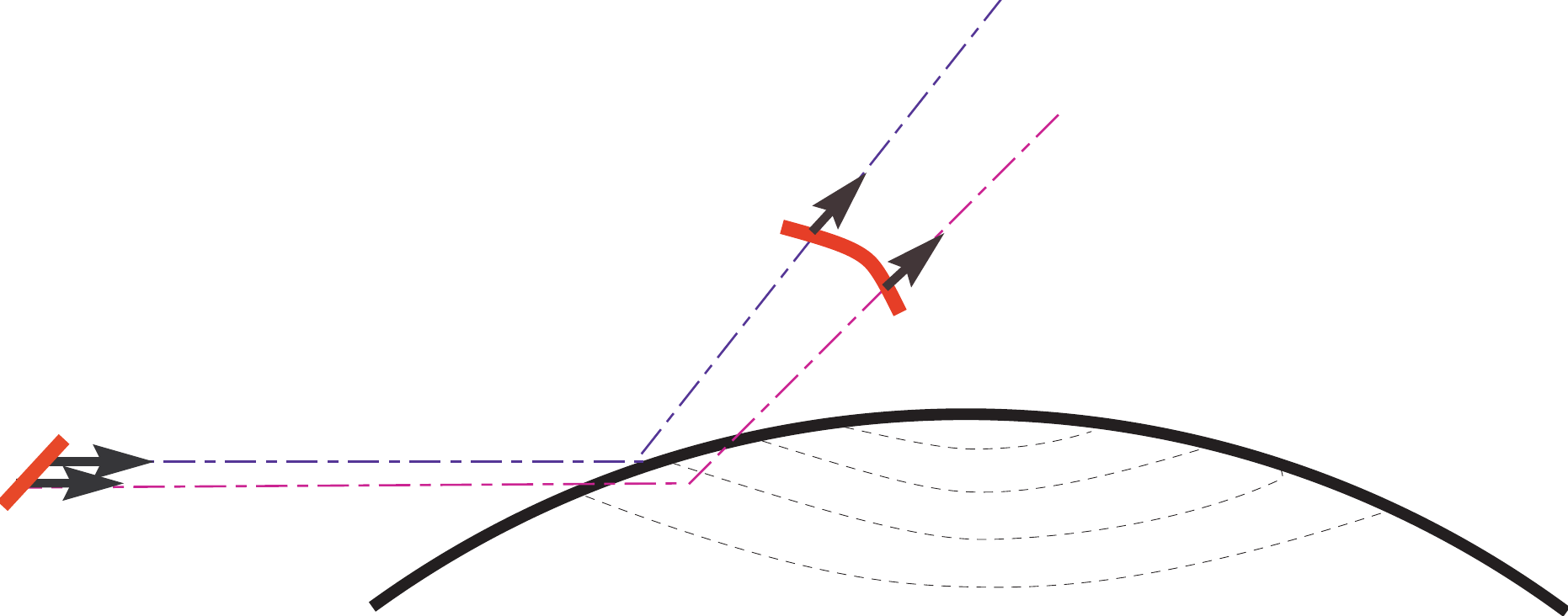}
\caption{Weak expansion}\label{fig:weak}
\end{figure}

The periodic {\it Lorentz process}, much interesting to physics, was
introduced by H.~Lorentz in 1905, see \cite{L05}. The periodic Lorentz process models the motion of an
electron among a periodic array of molecules in a metal and at the same time, as we will see, it also serves as a deterministic model of Brownian motion. Thus the  periodic {Lorentz process} is the natural $\mathbb{Z}^d$-cover of the
above-described toric billiard. More precisely: consider $\Pi:\mathbb{R}^d
\to \mathbb{T}^d$ the factorisation by $\mathbb{Z}^d$.	Its fundamental domain
$\itD$ is a cube (semi-open, semi-closed) in $\mathbb{R}^d$, so $\mathbb{R}^d =
\cup_{z \in \mathbb{Z}^d} (\itD+z)$, where $\itD + z$ is the translated fundamental
domain.
We also lift the scatterers to $\mathbb R^d$. Define the phase space of the Lorentz flow as $\tilde M = \tilde Q\times S_{d-1}$, where $\tilde Q = \cup_{z \in \mathbb{Z}^d} (Q+z)$ and that of the Lorentz map by $\partial \tilde M = \partial \tilde Q \times S_{d-1}^+$.  In the non-compact space $\tilde M$ the dynamics is denoted by  $\tilde S^t$  and the billiard map on $\partial \tilde M$  by	$\tilde T$. Their natural invariant measures on $\tilde M$ and $\partial \tilde M$ are $\tilde \mu$ and $\tilde \mu^\partial$, respectively (N. B.: they are infinite measures). The periodic Lorentz  process is the trace in the configuration space of the previous dynamical system, or in other words it is the natural projection $L(t) = L(t;x), t \in \IR_+$  to the configuration space $\tilde Q$ (or $\{L^\partial_{n}| n \in \ZZ_+$\}).


Here we will restrict our discussion to periodic Lorentz processes. We note that people also often discuss its non-periodic variants like for instance the random Lorentz process with random configurations of scatterers.

The {\it free flight
  vector} $\tilde\psi: \tilde M \to \mathbb{R}^d$ is defined as follows:
$\tilde \psi(\tilde x)=\tilde q(\tilde T\tilde x)-\tilde q (\tilde x)$. $|\tilde\psi|: \tilde M \to \mathbb{R}$ is called {\it the free flight function} or just the free flight.

\begin{definition}
  A semi-dispersing billiard (or the corresponding Lorentz process) is said to have {\it finite horizon} (FH) if the free flight vector
  is bounded.  Otherwise the system is said to have {\it infinite horizon} ($\infty$H).
\end{definition}


We note finally that systems of hard balls interacting via elastic collisions are isomorphic to semi-dispersing billiards. These are dispersing only if the number of the balls is two (as to a survey on hard ball systems see \cite{Sz08}).

\subsection{Hierarchy of stochastic properties}\label{subsec:prop}
For a deterministic dynamical system with a fixed invariant measure one can formulate an -- almost -- hierarchy of stochastic properties, stochastic behaviours which one can try to establish. Without aspiring after completeness let us list some  main levels for them:
\begin{enumerate}
\item hyperbolicity (existence of invariant cone fields (cf. \cite{W85}) or  existence of -- stable and unstable -- invariant manifolds);
\item non-zero Lyapunov exponents;
\item positive entropy;
\item positivity of ergodic components;
\item ergodicity;
\item qualitative forms of mixing (including K-property);
\item Bernoulli property;
\item quantitative forms of mixing;
\item central limit theorems  (CLT) or other forms of weak limit theorems;
\item local versions of limit theorems.
\end{enumerate}
For simplicity, properties (1-7) will be called as qualitative ones and properties (8-10) as quantitative ones; physicists are particularly interested in quantitative properties since they can be compared with measurements.

For completeness, we mention further interesting features though the discussion of related results goes beyond the scope of this survey.
\begin{enumerate}[resume]
\item strong or almost sure versions of limit theorems;
\item large deviation results;
\item stochastic properties related to mass or energy transport (for instance to Fourier law).
\end{enumerate}
Finally, for billiard type models or for the Lorentz process there are other intriguing 'geometric' questions (lying actually outside the hierarchy) like
\begin{enumerate}[label=(\alph*)]
\item calculation of the mean free path;
\item description of the horizon structure.
\end{enumerate}
An important question is whether the above listed properties do or do not differ for the discrete time maps (the billiard map $T$) and for the time-continuous billiard flow $S^t$. In many cases they are analogous but some of them requires additional effort. For instance, as explained in Section 6.8 of \cite{ChM06}, already the ergodicity of the flow does not follow automatically from that of the map.
In this paper for simplicity of exposition we almost exclusively focus on results related to the billiard map.
As said before, if the scatterers of a billiard (or of  a Lorentz process) are convex (i.~ e. the pieces of the boundary $\partial Q$ are concave inward $Q$), then the billiard is called {\it semi-dispersing} and if they are strictly convex, then it is called {\it dispersing}. The discussion of this paper is focused on the theory of semi-dispersing or dispersing billiards.

\subsection{Theory of planar billiards: some key points}\label{subsec:planekey}
\subsubsection{Qualitative properties}
Sinai in his celebrated 1970 paper \cite{S70} created the theory of planar dispersing billiards and -- by also assuming the finiteness of horizon -- proved their ergodicity and Kolmogorov mixing property (properties up to (6)). For doing so he had to extend the theory of planar smooth hyperbolic dynamical systems to two-dimensional hyperbolic systems with singularities.
  In particular, he showed that local, i. e.~ smooth versions of stable and unstable invariant manifolds do exist a.~e.  though they can be arbitrarily small. His analysis got extended to the infinite horizon case in \cite{BS73}. Based upon these results Gallavotti and Ornstein could also prove in \cite{GO74} the Bernoulli property of these billiards (property (7); we note that the work \cite{ChH96} provides general results and beyond that a nice survey of earlier ones).

\subsubsection{Quantitative properties}
The next major step was the proof of quantitative mixing properties and consequently that of the CLT through the construction of suitable Markov approximations  -- for FH models, at least (properties (8-9)).
\footnote{It is worth noting that, by having an appropriate correlation decay bound, the proof of a CLT is more or less standard though establishing the non-singularity of the limiting covariance of the diffusive limit may require additional  arguments, cf. \cite{B00}.}
 First Bunimovich and Sinai completed  this task by designing Markov partitions for a class of FH dispersing billiards, cf. \cite{BS80,BS86,BS81}. Ten years later they -- jointly with Chernov, \cite{BChS90,BChS91} -- relaxed the concept of Markov partition by introducing approximations via so-called Markov-sieves and, moreover, they also extended the previous results to a wider class of planar dispersing billiards. An important claim in these results was that for piecewise H\"older functions on the phase space one has stretched exponential decay of correlations. Then Young -- in \cite{Y98} -- could construct Markov extensions by her tower method. The advantage of her approach was not only that she could also obtain the optimal exponential decay of correlations for piecewise H\"older functions. Beyond that  she introduced a set of axioms that represented a framework for multidimensional models as well and furthermore covered or provided an access to dissipative systems -- like logistic maps (whose hyperbolicity was exploited earlier in e. g.~ \cite{J81,BC85}), the H\'enon map (whose hyperbolicity was exploited earlier in e.~g. \cite{BC91}), and the Lorenz map  as well.

Chernov, in \cite{Ch99}, extended the results of \cite{BChS90,BChS91}  to $\infty$H planar dispersing billiards for piecewise H\"older functions. The free jump vector in this case is, however, not H\"older and therefore his analysis did not provide yet a limit law for the corresponding Lorentz process. This was finally achieved in the following works.
The author, jointly with Varj\'u, could also establish local versions of the CLT for piecewise H\"older functions (property (10)). The FH case is treated in \cite{SzV04} and the $\infty$H case in \cite{SzV07}; the latter result also implied the verification of an earlier conjecture of Bleher \cite{B92} claiming that for the free flight vector a Gaussian limit theorem holds though with a non-standard scaling $\sqrt {n \log n}$ rather than with the standard $\sqrt n$ one. Finally for the free flight vector Chernov and Dolgopyat \cite{ChD09B}, could also derive exponential decay of correlations.

After this brief survey of results we recommend the interested reader to look into the excellent monograph \cite{ChM06} for a rather complete treatment of planar chaotic billiards.

An essential question, where, for instance, the behavior of the billiard map did -- until 2007 --  not imply that of the flow, was the bound on the correlation decay rate of piecewise H\"older functions. As said above, \cite{Y98} established the optimal decay rate for the map. It, however, did not imply an analogous result for the flow. It was only in 2007 that Chernov \cite{Ch07} could obtain a stretched exponential bound for the planar FH  dispersing billiard flow with smooth scatterers, and in 2015 that Baladi, Demers and Liverani \cite{BDL15} were able to achieve an exponential bound.

\subsubsection{Some examples of interaction of physics and mathematics}\label{subsub:inter}
\begin{enumerate}
\item In 1979 Bunimovich found an unanticipated phenomenon that since then has been arousing the interest of both mathematicians and physicists (theoretical, experimental and computational ones alike) demonstrating that chaotic behaviour is possible in a billiard which is nowhere dispersing: the Bunimovich stadium. \cite{ChM06} provides a detailed background on its mathematics so here I only add two remarks:
    \begin{enumerate}[label=(\roman*)]
\item  \cite{W86} gave a new insight into this phenomenon by also widely extending the class of billiards providing this behaviour;
\item \cite{B90} discovered that the local ergodicity theorem formulated for semi-dispersing billiards can even be used to settle the ergodicity of these nowhere dispersing billiards. Consequently many important results for semi-dispersing billiards can be transmitted to these new  types of billiards.
\end{enumerate}
\item The successes of mathematics described in the previous point had led to research on more sophisticated models, like Lorentz processes with stochastic reservoirs at the boundaries, the Gaussian iso-kinetic Lorentz process, Lorentz process in an external field, etc. Then Barra and Gilbert, \cite{BG07} suggested a billiard with a geometric bias and argued that it sustains a steady current. Motivated by \cite{BG07}  Chernov and Dolgopyat, \cite{ChD08} gave a rigorous proof and parallelly they were also suggesting that -- added to the steady drift -- there may also arise log-periodic oscillations. And in \cite{BChG07} it turned out that they indeed appear.
\item Similar example is  the recurrence of the Galton board model studied by Chernov and Dolgopyat \cite{ChD09C}, a surprising phenomenon that got obtained and explained on the basis of rigorous study.
    \item In general, the method of averaging, so successful in celestial mechanics and, in general, in investigations of non-linear dynamical systems, has become rewarding in billiard models, too. A spectacular example is the work \cite{ChD09} of Chernov and Dolgopyat on a mechanical -- actually a billiard -- model of Brownian motion (cf. Subsection \ref{subsec:brown}).
        \item The wind-tree model of the Ehrenfests' studied by Delecroix, Hubert and Leli\`evre, \cite{DHL14} is a billiard model and is not hyperbolic at all. Simply saying it is a planar Lorentz process with infinite horizon and periodically situated rectangular scatterers. The authors prove rigorously an astonishing result that was not foreseen by physical arguments, neither whose heuristic explanation would be easy:  independently of the size of the rectangular scatterers -- for typical angles of the initial velocity -- the polynomial diffusion rate of the system is $2/3$ (meaning, roughly speaking, that typical trajectories occasionally go from the origin as far as $T^{2/3}$ where $T$ denotes time).
\end{enumerate}

\subsection{Theory of multidimensional billiards: key points}\label{subsec:multikey}
\begin{figure}[h]
\includegraphics[scale=0.4]{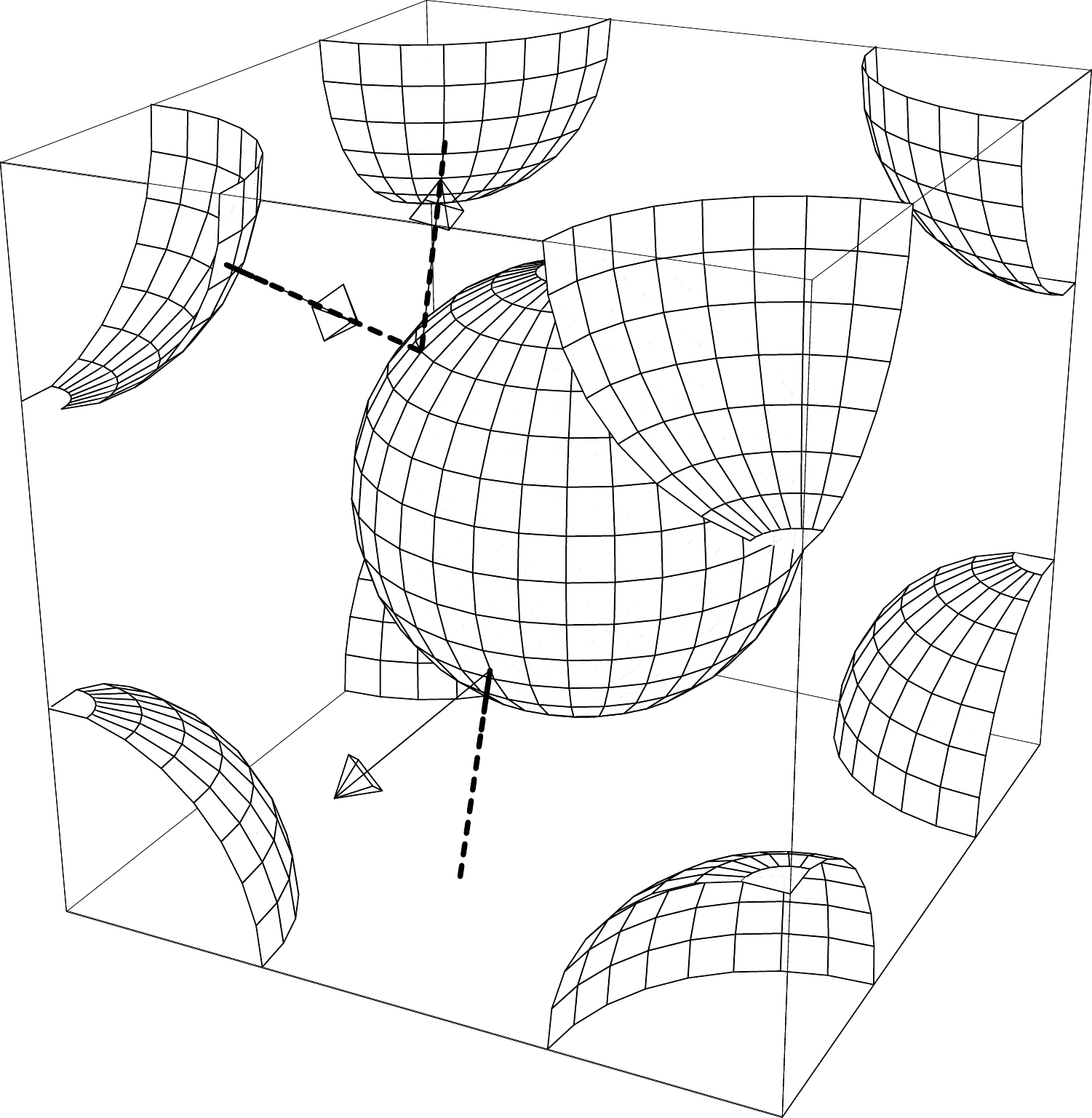}
\caption{3-dimensional dispersing billiard}\label{fig:3Ddispersing}
\end{figure}

Since the main part of this paper is devoted to multidimensional billiards, here we are satisfied by giving a snapshot of the main results. Further details will be provided in Sections \ref{sec:multiqual} and \ref{sec:multiquan}. Figure \ref{fig:3Ddispersing} shows a 3-dimensional dispersing billiard configuration.

\subsubsection{Qualitative properties}
\begin{enumerate}
\item Chernov \cite{Ch83} initiated the theory by showing that -- for a semi-dispersing billiard satisfying some mild conditions -- local stable and unstable invariant manifolds exist for a.~e. phase point (warning: the dimensions of the invariant manifolds may not be maximal, i.~e.  may be less than $d-1$). This fact also implied the positivity of the metric entropy. We note that his analysis was using the evolution equations for convex fronts elaborated in \cite{S79,ChS82}.
\item The next much important development was the work \cite{ChS87} treating semi-dispersing billiards, in general. Its authors introduced the notion of a sufficient phase point, i.~e.  one which -- in its collision history -- encounters all possible degrees of freedom of the system. The authors showed that -- under some conditions also including the novel Chernov-Sinai Ansatz -- an open (!) neighborhood of a sufficient phase point belongs to one ergodic component (``a local ergodicity theorem").
\item The twin papers \cite{BChSzT02} and \cite{BChSzT03} discovered a pathology appearing in case of orbits suffering at least two tangential collisions under particular angles; an example  is shown on Figure \ref{fig:pathol}. This phenomenon actually  implies that images of singularities are not smooth submanifolds contrary to earlier conviction that they are. These works could also establish that the theory of \cite{ChS87,KSSz90} is nevertheless valid for algebraic billiards, i.~e.  for billiards with piecewise algebraic boundaries. Moreover, \cite{BChSzT03}  provided a self-contained geometric description of unstable manifolds together with proofs of fundamental regularity properties of multidimensional dispersing billiards.

\begin{figure}[hbt]
\centering
\epsfxsize=10cm
\epsfysize=6cm
\epsfbox{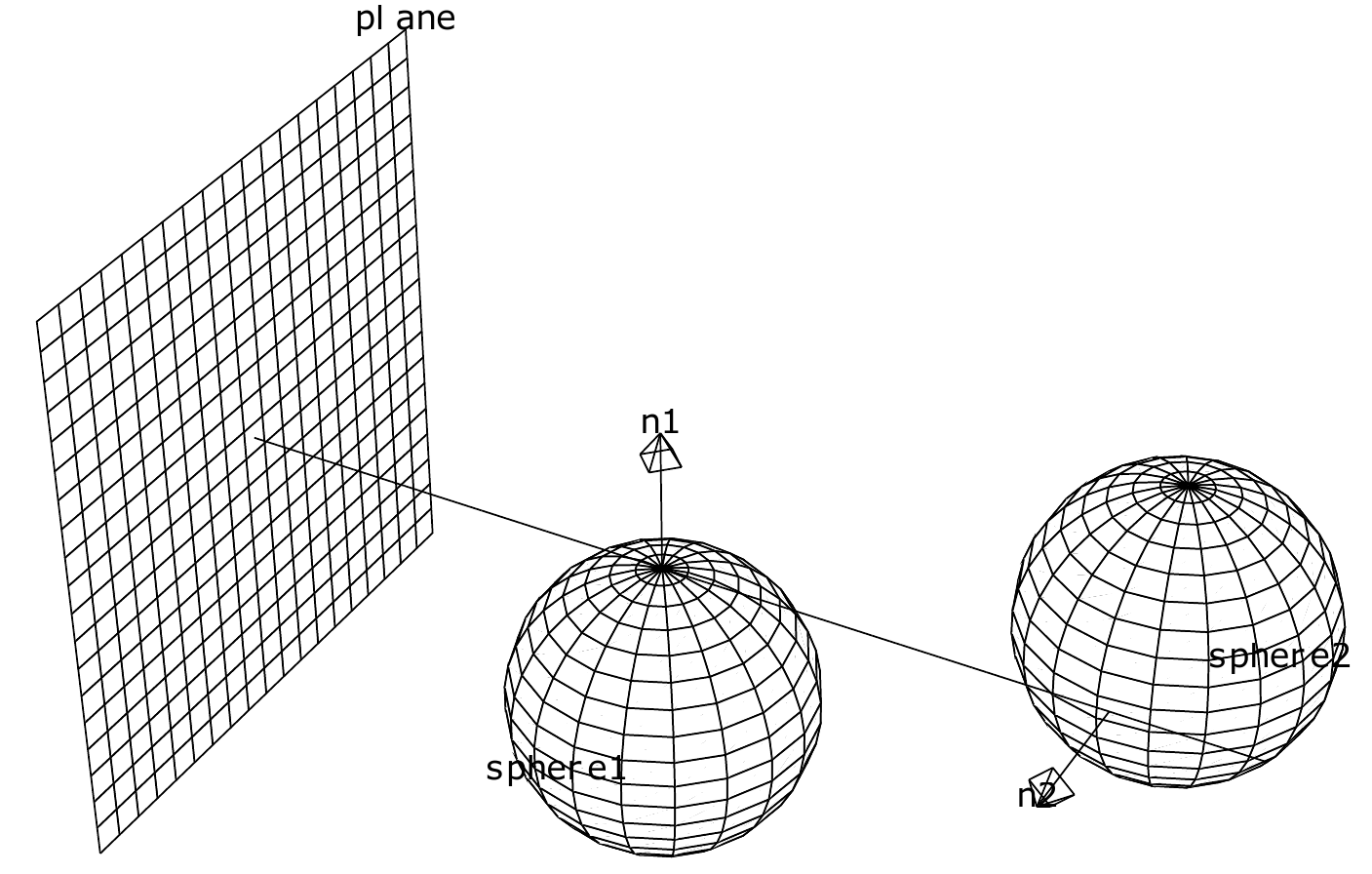}\caption{The studied billiard configuration}\label{fig:pathol}
\end{figure}

\item The local ergodicity theorem of  \cite{ChS87} also opened the way to establishing a major achievement of the theory: the celebrated Boltzmann-Sinai ergodic hypothesis claiming that systems of elastic hard balls are ergodic (immediately implying K-mixing, too) modulo the trivial invariants of motion.
     (The first step was reached in \cite{KSSz91}, the last sensational one through a series of works by Sim\'anyi concluding in \cite{S13}; as to a recent survey of the main ideas and problems see \cite{Sz08}.) At this point it is worth noting that hard ball systems are isomorphic to algebraic billiards.
    \end{enumerate}
 \subsubsection{Quantitative properties}
 \begin{enumerate}[resume]
\item The main results of Young's seminal work \cite{Y98} offered a framework suitable for the study of higher dimensional hyperbolic dynamics and while doing so they also provided optimal correlation decay rate for some centrally interesting one- or two-dimensional hyperbolic systems. At this point we remind the reader of Sinai's classical billiard philosophy: for billiards the classical tools of smooth hyperbolic dynamical systems work since {\it expansion prevails partitioning} (the latter caused by the singularities). Various quantitative formulations of this idea had been appearing from the start of the theory and they led in \cite{Y98} to explicit technical {\it growth lemmas}, a tool that later became fundamental (cf.  \cite{Ch99,Ch99b,BBT08,BT08} and to several further works of Chernov with coauthors, mainly with Dolgopyat).
       As a result, in principle, in multidimensional billiards it is a paramount task to control the sizes and the shapes of local stable and unstable manifolds. This framework of \cite{Y98} helped Chernov \cite{Ch99b} to treat this problem for arbitrary dimensional smooth hyperbolic systems with singularities and thus for systems where -- unlike in billiards -- the derivative map is bounded.

\end{enumerate}
Similarly to the properties listed in subsection \ref{subsec:prop}, there also arise intriguing `geometric' properties in the multidimensional case.
\begin{enumerate}[resume]
\item Chernov \cite{Ch97,Ch00} gave formulas for the entropy of hyperbolic billiard maps and flows (the second work also provided a nice survey of the topic). He also calculated the mean free path providing the ratio of the two entropy values.
\item Dettmann \cite{D14} described the rich world of $\infty H$ multidimensional billiards and formulated three nice conjectures for the asymptotic tail probabilities of the free flight; two geometrical ones got established in  \cite{NSzV15} while the third dynamical one is still open.
\end{enumerate}

\section{Qualitative properties of multidimensional billiards}\label{sec:multiqual}
We start by summarizing some basic notions from \cite{BChSzT03} and for more details we also refer to that paper.
Consider a billiard on $\IT^d$ and -- for simplicity of exposition -- {\it assume that it is dispersing with smooth scatterers and with finite horizon}. First we are going to introduce the notions of fronts, convex fronts, $u$-manifolds and invariant manifolds for the billiard map $T$. They can also be introduced for the flow but here we will not need this generality.
\subsection{Fronts}
At this point it is more natural to consider the phase space $M$ of the flow.
Take a smooth 1-codimensional submanifold $E \subset Q$ and add the unit normal vector $v(r)$ of this
submanifold at every point $r \in E$ as a velocity, continuously. (Consequently, at every point
the velocity points to the same side of the submanifold $E$.) Then we say that
$$
{\itW}=\{(q,v(q))\vert q\in E\}\subset {M},
$$
is a wave front or simply a {\it front} if $v:E\to \sd$ is continuous (smooth) and $v\perp E$ at every point
of $E$.
The derivative of this function $v$, denoted by $B$, plays a crucial role:
$dv=Bdq$ for tangent vectors $(dq,dv)$ {of the front} ${\itW}$. $B$
acts on the tangent plane ${\itT}_qE$ of $E$, and takes its values
from
the tangent plane ${\itJ}={\itT}_{v(q)}S_{d-1}$ of the velocity
sphere.
These are both naturally embedded into the configuration space $Q$ and
can be identified through this embedding. So we just write
$B:{\itJ}\to{\itJ}$ with $B$ being the curvature operator of
the submanifold $E$.
Obviously, $B$ is symmetric.
\begin{remark}\label{rem:dfront}
It is important to observe that fronts defined in the phase space as above can be identified by their traces $E$ in the configuration space modulo their orientation. In what follows it will be convenient to assume that $B$ is non-negative and then call the front a convex one (\cite{ChM06} calls those with $B > 0$ dispersing fronts).
\end{remark}

One can also introduce {\it fronts for the map} as their traces in $\partial M$ taken in the following way: for every point of $(q,v) \in \itW$ of the convex front $\itW$  one takes its evolution under the flow exactly until the next collision, and then applies the collision law, i.~e. the collision operator to this incoming phase point to obtain the corresponding outgoing phase point, a point of the front of the map.

\subsection{Singularities}\label{subsec:sing}

The billiard map $T$ is discontinuous
at pre-images of tangential reflections. Indeed, consider the set
of tangential reflections:
$$
{\itS}_0: = \{ (q,v)\  | \ \la v,n(q) \ra =0 \}
$$
(i.~e.  the boundary $\partial(\partial M)$ of the phase space). Its
pre-images are:
$$
{\itS}_k= T^{-k} {\itS}_0 \quad (k>0).
$$
For later use it is useful to denote
${\itS}_{[k]}=
\cup_{i=1}^{k} {\itS}_i$.
 The map $T$ is discontinuous precisely at the
points of ${\itS}_1(={\itS}_{[1]})$. Similarly the billiard flow is discontinuous whenever the trajectory hits $\itS_0$.

When hitting a singularity, i.~e.  in points of ${\itS}_1$,  the derivative $dT$ explodes. Therefore, for obtaining appropriate regularity estimates, for instance distorsion bounds, fundamental in the theory of hyperbolic dynamical systems, one introduces additional, so-called secondary singularities by further partitioning the phase space into so-called homogeneity layers:
\begin{align*}
    I_k &=& \{(q,v)\in \partial M |\, &(k+1)^{-2}< \la v,n(q)\ra <k^{-2}\} \quad
     \nn \\
    I_0 &=& \{(q,v)\in \partial M |\, &\la v,n(q)\ra > k_0^{-2}\}
\end{align*}
Here the integer $k_0$ is large enough.
Thus new boundary pieces appear in the phase space $\partial M $, denote their union as
$$
\Gamma_0 = \cup_{k=k_0}^{\infty}  \{(q,v) | \la v,n \ra =
k^{-2} \}
$$
Correspondingly, the countably many manifolds in the set $\Gamma_{1} = T^{-1}
\Gamma_0$ are the so called {\it secondary singularities}.

For a higher iterate of the dynamics, $T^n$, the primary and secondary
singularities are, respectively:
$$
{\itS}_{[n]}= {\itS}_{1}\cup T^{-1}{\itS}_{1}\cup\cdots T^{-n+1}{\itS}_{1}; \quad
\Gamma_{[n]}= \Gamma_{1}\cup T^{-1}\Gamma_{1}\cup\cdots T^{-n+1}\Gamma_{1}.
$$

\subsection{$u$- and $s$-manifolds, invariant manifolds}

A front will be called {\it convex/diverging} whenever $B^+$ is positive
semi-definite ($B^+ \ge 0$). The {\it convex cone} consists of those tangent vectors $\delta x$
that are tangent to some convex front. A natural choice for the unstable cone field $\{C^u_x | x \in M\}$, invariant under the map $T$ (cf. \cite{W85}) is the family of convex cones still sometimes other choices are preferred. Though convex fronts remain convex under time evolution, they may get cut by singularities.

We say that a $d-1$-dimensional smooth manifold $\itW \subset Q$ (cf. Remark \ref{rem:dfront}) is a $u$-{\it manifold} if all its tangent planes lie in the (appropriately defined) unstable cone. A basic fact is that a u-manifold $\itW$  is an {\it unstable invariant manifold} or briefly an {\it unstable manifold} if each map $\{T^n | n \le -1\}$ is smooth on $\itW$ -- by taking into account the secondary singularities, too (sometimes such unstable manifolds are called (strongly) homogeneous unstable manifolds)

The a.~e.  existence of invariant manifolds for multidimensional semi-dispersing billiards was first established in \cite{Ch83}. For dispersing billiards their existence follows simply from the works \cite{P92,Y98}. In fact, since Sinai's original work \cite{S70}, the arguments behind this existence are based on a simple application of Borel-Cantelli lemma: -- by cutting out convergently decaying neighborhoods of higher and higher order singularities $\itS_n$ of $T^{-1}$ -- $u$-manifolds through phase points lying in the remaining domain have controllable neighborhoods not meeting any sufficiently high order singularity from $\itS_{[n]} \cup \Gamma_{[n]}$. In this small neighborhood then one can refer to the theory of smooth hyperbolic systems. It is worth stressing that these results do not require strong assumptions on the smoothness of singularities.

\subsection{Local ergodicity}
{\it Local ergodicity} at a point $x \in \partial M$ means that an open neighborhood of $x$ belongs to one ergodic component. The importance of this concept can be explained as follows. For a smooth dynamical system with non-zero Lyapunov exponents Pesin's theory \cite{BP07} ensures that the ergodic components of the system have positive measure (and hence there is at most a countable number of them). If one wishes to establish ergodicity, then the mere information  on the positivity of the components is not sufficient. For connecting the ergodic components in a hyperbolic billiard -- as this was initiated in \cite{S70,BS73,ChS87} -- information on their topology (e.~g.  openness)  and on the abundance of phase points possessing ergodic neighborhoods  is also needed. Going ahead in our exposition we remark that, if every phase point $x$ had an open neighborhood of $x$ belonging to one ergodic component, then the ergodicity of the system would be implied automatically. In dispersing billiards the case is almost so nice since the complement of the set of `good' points (and points in $\partial M^*$ (cf. \ref{equ:phase})  are certainly good by Theorem \ref{thm:leth}) can not separate possible ergodic components.

Let us briefly mention some key results of the theory of the local ergodicity theorem in the multidimensional case.
 \begin{itemize}
\item The first theorem on local ergodicity was formulated in \cite{ChS87} for semi-dispersing billiards.
\item Its generalization, called the transversal fundamental theorem, appeared in \cite{KSSz90}. At some points this work also contributed to the simplification and clarification of the original arguments.
\item \cite{LW95} treated local ergodicity for piecewise smooth, hyperbolic Hamiltonian systems. Thus it does not cover billiards where -- as we know -- the derivatives of the map explode. Nevertheless an important merit of the work was that it embedded billiard theory into the symplectic setup. It is worth mentioning that in these systems there is no strong anisotropy and thus images of the singularities are still smooth.
\item As said above, the twin papers \cite{BChSzT02} and \cite{BChSzT03} discovered that -- as a consequence of the strong anisotropy of billiard collisions -- the theorems of \cite{ChS87,KSSz90} are so far restricted to algebraic billiards, i.~e.  to billiards with piecewise algebraic boundaries only. (Figures \ref{fig:strong} and \ref{fig:weak} illustrate the {\it strong expansion} in the plane determined by the incoming velocity and the normal vector of the scatterer and the {\it weak expansion} in directions orthogonal to it, respectively.)
   \end{itemize}

Let us turn to a technical formulation of the local ergodicity theorem. First we introduce some more notation. For any $n \in \IN$,  $\Delta_n$ stands
for
the set of doubly singular phase points up to order $n$, i.e. $x\in\partial M$
belongs to $\Delta_n$ whenever there are indices $k_1\ne k_2$, $|k_i| \le n$
such that both $T^{k_1}x$ and $T^{k_2}x$ are elements of $\itS_0$. We will need the following sets:
\bea\label{equ:phase}
\partial M^0:= & & \partial M \setminus \bigcup_{n\in \ZZ}\itS_n \nn \\
\partial M^*:= & & \partial M \setminus \bigcup_{n=1}^{\infty}\Delta_n \nn
\eea
$\partial M^0$ is the set of phase points that never hit any singularity and $\partial M^*$ of those that hit at most one singularity during their whole evolution.
For dispersing billiards the conditions of the local ergodicity theorem are relatively simple. The reader is reminded that for simplicity we assume that the boundaries of the scatterers are smooth, i.~e.  there are no corner points (they may cause additional problems, in particular in obtaining appropriate complexity bounds, cf. \cite{dST14}).

\begin{theorem}\label{thm:leth}
Assume a dispersing billiard is algebraic, i.~e.  the boundaries of the scatterers are algebraic. Then any $x \in \partial M^*$ has an open neighborhood which belongs to one ergodic component.
\end{theorem}

\begin{corollary}
Under the conditions of the theorem, any connected component of the configuration space of a dispersing billiard defines an ergodic component of the billiard system.
\end{corollary}

\begin{remark}
\begin{enumerate}
\item The condition on algebraicity of $\partial Q$ follows from a weaker technical one: there exists a constant $L \in \IR_+$such that for every $N \in \IN$ the singularity set $\cup_{n \le N} \itS_n$ is Lipschitz decomposable with constant $L$ (i.~e.  it is a finite collection of graphs of Lipschitz functions, cf. Remark 3.6 in \cite{BChSzT02}).
\item As said above, we restrict ourselves to dispersing billiards therefore the Chernov-Sinai Ansatz (Condition 3.1 from \cite{KSSz90}) holds automatically. (For semi-dispersing billiards, for instance like hard ball systems, its verification is far from trivial, cf. \cite{S13}.)
\item The local ergodicity theorem follows from the so-called fundamental theorem whose formulation is definitely more involved therefore we omit it and refer to the papers listed above.
\item For experts it was a surprise that the algebraicity of scatterers appeared as a condition of the local ergodicity theorem and there is a general agreement that it is not necessary. An encouraging fact is the result of \cite{BBT08} to be discussed in the next section where only $C^3$-smoothness of the boundaries of the scatterers is assumed.
\end{enumerate}
\end{remark}

\begin{problem}\label{nonalgebraic}
The weakening of the conditions of the local ergodicity theorem (i.~e. those on algebraicity or Lipschitz decomposability of the boundary pieces of the scatterers) is a major mathematical question of the theory.
\end{problem}
\noindent I say it is mathematical since the examples most interesting from the point of view of physics: hard ball systems and their variants are all algebraic. It is worth mentioning that examples (see \cite{BR98}) of multidimensional nowhere dispersing billiards are also algebraic billiards (cf. Point 1 of Subsection 2.3.3).

\subsection{Complexity bound}
As we have already mentioned, Sinai's basic billiard philosophy was: {\it expansion prevails partitioning} and this underlying principle has been used throughout. As an explicit condition it appeared when proving a quantitative property: exponential correlation decay in two-dimensional systems in \cite{Y98}, cf. hypothesis (H4) in subsection 7.1. Afterwards in \cite{Ch99,BT08,D09} the complexity assumption got formulated in the {\it multidimensional} setup and -- interestingly enough -- was also used in obtaining local ergodicity, a qualitative property.

Following \cite{BT12} we introduce the $n$-step complexity of the billiard map $T$. Consider the connected components $\{O_j| 1 \le j \le K\}$ of the set $\partial M \setminus \itS_{[1]}$.  The map $T$ itself is, of course, continuous on each of these open components. Next, for every $n \ge 1$ and for each of the multi-indices ${\bf i} = (i_1, \dots , i_n) \in \{1,  \dots, K\}^n$ let $T^n_{\bf i} =: T^{i_n} \circ \dots \circ T^{i_1}$ be defined on the set $O_{\bf i}$ where  the subsets $O_{\bf i}$ are defined inductively as follows: $O_{(i_1)} = O_{i_1}$ and
$$
O_{\bf i} = O_{(i_1, \dots, \i_n)} = \{ x \in O_{i_1}\  {\rm such\ that}\  T_{i_1}x \in O_{(i_2, \dots, i_n)}\}.
$$
\begin{definition}
For every $x \in \partial M$, the $n$-step complexity index $m_n(x)$ of $x$ is the number of the possible distinct multi-indices ${\bf i} \in \{1,  \dots, K\}^n$ such that $x \in {\it Cl}\ (O_{\bf i}) $. Then the $n$-step complexity $m_n$ of the singularity set of $T$ is the supremum of $m_n(x)$ over $x \in \partial M$.
\end{definition}
 Denote, moreover, the minimal expansion rate along unstable vectors by $\Lambda_{\rm min}$.
\begin{condition}[Sub-expansion complexity condition]
 The billiard satisfies the sub-expansion complexity condition if for some $1 < \lambda < \Lambda_{\rm min}$ one has  $m_n = \cO(\lambda^n)$.
\end{condition}

A stronger form of the complexity assumption is the following one:
\begin{condition}[Sub-exponential complexity condition]
 The billiard satisfies a sub-exponential complexity condition if for every $1 < \lambda $ one has  $m_n = \cO(\lambda^n)$.
\end{condition}

Translating the sub-expansion complexity condition to words just says the following: our basic tools are unstable manifolds and they are expanded during the first $n$ iterations by a factor of at least $\Lambda^n$; during this time the images of an unstable manifold are cut by singularities; the condition requires that the number of pieces be less than the minimal expansion during these steps. Precise technical formulations of this intuition were given at several places, for instance in \cite{Ch99,BT08}. In \cite{D09} the growth lemma is directly derived from the sub-expansion complexity bound.

As noticed in relation to Problem \ref{nonalgebraic} -- for multidimensional uniformly hyperbolic systems with singularities -- Bachurin, B\'alint and T\'oth, \cite{BBT08} could derive a  fundamental theorem by assuming proper alignment, appropriate distortion bounds and by proving the growth property.

\section{Quantitative properties of multidimensional billiards}\label{sec:multiquan}
\subsection{Exponential correlation decay and CLT}
Unfortunately there is only one -- actually  a most remarkable -- result by B\'alint and T\'oth \cite{BT08} and, moreover, it is conditional: it presupposes a sub-expansion complexity bound.
\begin{theorem}[\cite{BT08}]
Consider a finite horizon dispersing billiard on $\IT^d$ with $C^3$-smooth scatterers and assume it satisfies the sub-expansion complexity condition. Then the dynamics enjoys exponential correlation decay and -- for non-coboundary functions -- satisfies the central limit theorem with a non-singular limiting covariance matrix.
\end{theorem}

Though this theorem is conditional, nevertheless it represents a much significant step in the quantitative study of multidimensional billiards. Therefore it is a big challenge of the theory to clarify the situation around the complexity assumption much the more since in physical questions the control of correlation decay is crucial.
\begin{problem}
Construct  a   finite horizon dispersing billiard on $\IT^d: d \ge 3$ with $C^3$-smooth scatterers that satisfies the sub-expansion complexity condition.
\end{problem}
(I note that it is easily possible that the shapes of the scatterers in such an example  will not be spherical.)

Non-constructive forms of the problem are the following ones:
\begin{problem}
\begin{enumerate}
\item Show the existence of a   finite horizon dispersing billiard on $\IT^d: d \ge 3$ with $C^3$-smooth scatterers satisfying the sub-expansion complexity condition.
\item Prove that the set of $d$-dimensional ($d \ge 3$) finite horizon dispersing billiards on $\IT^d: d \ge 3$ with $C^3$-smooth scatterers on  satisfying the sub-expansion complexity condition is generic -- in a reasonable sense of genericity.
\end{enumerate}
\end{problem}

\begin{figure}[h]
\includegraphics[scale=0.6]{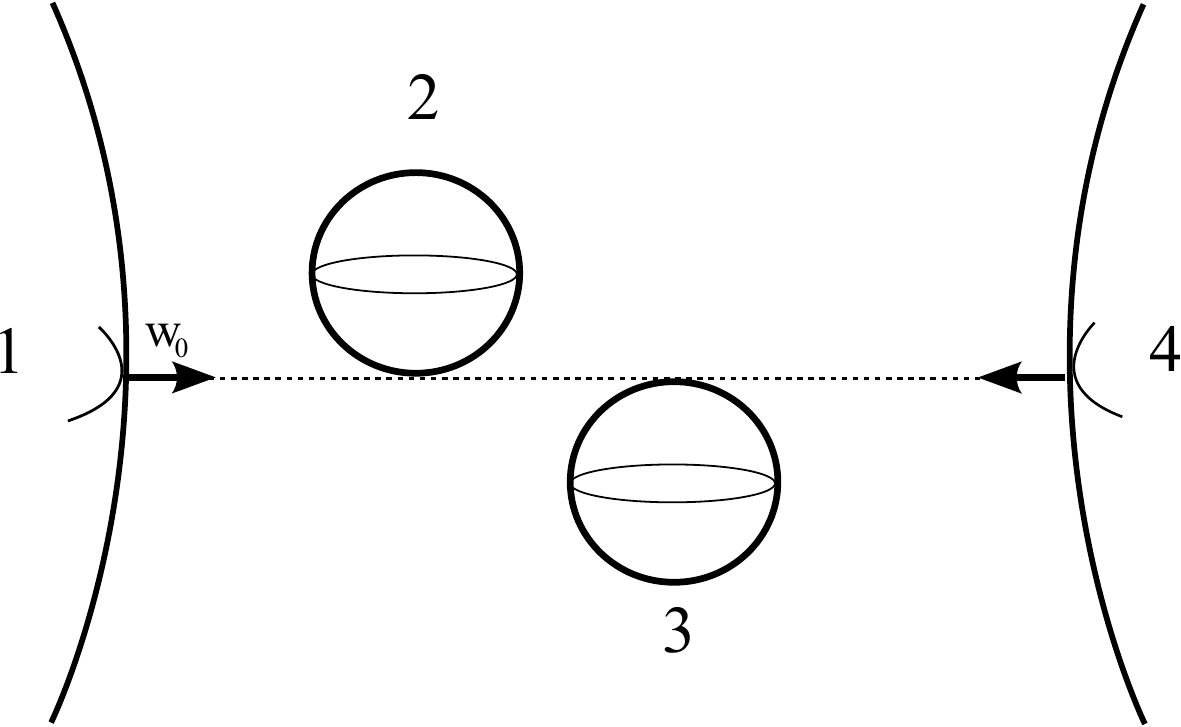}
\caption{A 3-D billiard configuration with exponential complexity}\label{fig:counterexample}
\end{figure}

From the other side there is an interesting example \cite{BT12}  -- already in dimension 3 -- of a finite horizon dispersing billiard where not only the complexity grows exponentially but its rate is larger than the minimal expansion rate implying that the sub-expansion complexity assumption does not hold for this billiard (see Figure \ref{fig:counterexample}). The existence of such an example does not contradict to exponential correlation decay and it is not surprising either as the following very rough intuition suggest this: the number of connected pieces  of the images of an unstable manifold as they are cut by singularities is approximately proportional to its average increase rate in the directions 'transversal to the singularities'. Indeed the fact that for planar dispersing billiards the complexity bound is linear is in accordance with the previous picture: in planar billiards the expanding direction is one-dimensional and thus in the `transversal direction' the increase rate is zero. Therefore the sub-expansion complexity condition, in general, certainly seems to be too strong. The work \cite{BT12} contains a detailed analysis of the complexity question and is highly recommended for further details.

As to the perspective for settling the uncertainty around the complexity issue an encouraging fact is that various forms of complexity assumptions also arise when studying piecewise expanding maps (cf. \cite{S00,T01,B01,B09}) and the intuition gained there may help here, too.
\subsection{Systems close to be planar}\label{subsec:brown}
Chernov and Dolgopyat \cite{ChD09} introduced a much interesting mechanical model for the study of Brownian motion and obtained fascinating results for it. Take a finite horizon planar dispersing billiard in a container  $Q \subset \IT^2$ with $C^3$-smooth scatterers and put two disjoint particles in $Q$. One of them: $(\cQ,\cV)$ is a hard disk of radius $r > 0$ and of mass $M \gg 1$ and the other one: $(\mathcal q,\mathcal v)$  is a point particle of mass $m=1$. The two particles move according to the laws of mechanics: elastic collision of the particles and specular reflection at the scatterers. Consequently the kinetic energy of the system is conserved. We can assume
\[
||\mathcal v||^2 + M||\cV||^2 = 1
\]
We note that this system is, in general, hyperbolic and ergodic (see \cite{SW89}). \cite{ChD09} is, however, treating the rescaled velocity process $\cV(t)$ of the heavy particle when $M \to \infty$. An important fact is that the time scale the authors are  studying is short compared to the mixing time of the two-particle system. In their limiting case $M \to \infty$ the heavy particle moves very slowly and thus it is almost still. Consequently the dynamics of the light particle is almost a planar dispersing billiard whose scatterers  are the original ones plus the immobile heavy particle. Technically speaking the system of the two particles together is isomorphic to a 4-dimensional semi-dispersing billiard. For it the phase space of the map is 6-dimensional. The technical advantage of the $M \to \infty$ limit is that then the unstable invariant cones of the two-particle system are asymptotically close to the unstable invariant cones of the aforementioned planar billiard.

\noindent {\bf Acknowledgements.} I am most grateful to the organizers of {\it The Dynamical Systems, Ergodic Theory, and Probability Conference Dedicated to the Memory of Nikolai Chernov}. My special thanks are due to P\'eter B\'alint and Imre T\'oth for their most careful reading of the manuscript and a bundle of substantial remarks and, moreover, for the figures. Thanks are due to the referee and to P\'eter N\'andori for their precious remarks. I am also indebted to Carl Dettmann and Thomas Gilbert for examples treated in Subsection \ref{subsub:inter}. This research was supported by Hungarian National Foundation for Scientific
Research grants No. K 71693, K 104745 and OMAA-92öu6 project.
\vskip5mm
\appendix
\centerline{\large APPENDIX}
The goal here is to illustrate the motto of the paper by extracting the most significant works of Kolya Chernov in creating and advancing the theory of multidimensional billiards.
\begin{itemize}
\item Following \cite{K79}, \cite{ChS82} works out the collision equations for multidimensional billiards;
\item \cite{Ch83} is devoted to the construction of invariant foliations by stable and unstable manifolds and to showing their transversality for semi-dispersing billiards. A conclusion of the results is the positivity of Kolmogorov-Sinai entropy.
\item In \cite{ChS87} the authors find and prove the local ergodicity theorem that gave a possibility and strong impetus among others to proving ergodicity for hard ball systems.
\item \cite{Ch99b} is dedicated to developing tools to control the shape of invariant manifolds for piecewise smooth
hyperbolic systems in high dimensions (these systems do not contain billiards since here the derivative map is assumed to be bounded).
\item The previous work together with \cite{Ch99} provide a much suitable form of growth lemmas, basic ingredients for treating quantitative properties;
\item The papers \cite{BChSzT02,BChSzT03} discovered a pathology appearing multidimensional billiards and leading to the non-smoothness of images of singularities. They provided conditions for the validity of the local ergodicity theorem and also gave a detailed analysis of geometric properties of  multidimensional billiards;
\item \cite{ChD09} introduced the methods of averaging and of standard pairs into the study of billiards and by using it they  could give a dynamical model of Brownian motion in a multidimensional billiard which is close to be planar.
\end{itemize}

\end{document}